\documentclass[12pt]{article}
\usepackage{amssymb}
\usepackage{amsmath}

\newcommand{\lap}{\mbox{$\bigtriangleup$}}

\newcommand{\ra}{{\mbox{$\rightarrow$}}}
\newcommand{\be}{\begin{equation}}
\newcommand{\ee}{\end{equation}}
\newtheorem{mrem}{Remark}
\newtheorem{mthm}{Theorem}

\newtheorem{mlem}{Lemma}

\newtheorem{thm}{Theorem}[section]

\begin{document}

\title{A maximum principle on unbounded domains and a Liouville theorem for fractional \emph{p}-harmonic functions}

\author{Wenxiong Chen \thanks{ Corresponding author, Partially supported by the Simons Foundation Collaboration Grant for Mathematicians 245486.} \quad Leyun Wu  }

\date{\today}
\maketitle
\begin{abstract} In this paper, we establish the following Liouville theorem for fractional \emph{p}-harmonic functions.

{\em Assume that $u$ is a bounded solution of
$$(-\lap)^s_p u(x) = 0, \;\; x \in \mathbb{R}^n,$$
with $0<s<1$ and $p \geq 2$.

Then $u$ must be constant.}

A new idea is employed to prove this result, which is completely different from the previous ones
in deriving Liouville theorems.

For any given hyper-plane in $\mathbb{R}^n$, we show that $u$ is symmetric about the plane. To this end, we established a {\em maximum principle} for anti-symmetric functions on any half space.
We believe that this {\em maximum principle}, as well as the ideas in the proof, will become useful tools in studying a variety of problems involving nonlinear non-local operators.
\end{abstract}
\bigskip

\medskip

{\bf Math subject classification} 35J60

\section{Introduction}

In $n$-dimensional Euclidean space $\mathbb{R}^n$, the fractional $p$-Laplacian is a fully nonlinear pseudo-differential operator defined by
\begin{eqnarray*}
(-\Delta)_p^s u(x)&=&C_{n,sp} \lim_{\varepsilon\rightarrow 0} \int_{\mathbb{R}^n \backslash B_{\varepsilon}(x)}\frac{|u(x)-u(y)|^{p-2}(u(x)-u(y))}{|x-y|^{n+sp}}dy\\
&=&C_{n,sp} PV \int_{\mathbb{R}^n }\frac{|u(x)-u(y)|^{p-2}(u(x)-u(y))}{|x-y|^{n+sp}}dy,
\end{eqnarray*}
where $C_{n,sp}$ is a constant depending on $n$, $s$, and $p$, and $PV$ stands for the Cauchy principal value.

In order the above integral to converge, we require that
$$
u \in C_{loc}^{1,1}\cap {\cal L}_{sp}
$$
with
$$
{\cal L}_{sp} =\{u\in L_{loc}^{p-1} \mid \int_{\mathbb{R}^n }\frac{|1+u(x)|^{p-1}}{1+|x|^{n+sp}}dx<\infty\}.
$$

We say that $u$ is a fractional \emph{p}-harmonic function in $\mathbb{R}^n$ in the classical sense if
\begin{eqnarray}\label{p-laplacian}
(-\lap)_p^su(x)=0, \;\;\;\;\;  x \in \mathbb{R}^n.
\end{eqnarray}

One of our main results in this paper is

\begin{mthm}(Liouville theorem)\label{main}
Assume that $u \in C_{loc}^{1,1}\cap {\cal L}_{sp}(\mathbb{R}^n)$ is bounded, if
\begin{eqnarray}\label{p-laplacian-1}
(-\lap)_p^su(x)=0 \text{~in~} \mathbb{R}^n,
\end{eqnarray}
then
\begin{eqnarray}\label{constant}
u\equiv C.
\end{eqnarray}
\end{mthm}
\begin{mrem}\label{counterexample}
For unbounded functions, there are obvious counter examples. For instance, if $u(x)=x_i$, $i=1, 2, \cdots, n$, then one can easily verify that
$$
(-\Delta)_p^s u(x) =0, \;\;\; \forall \, x \in \mathbb{R}^n.
$$
\end{mrem}

The well-known classical Liouville's theorem for harmonic functions states that

{\em If $u$ is bounded from below or from above and
$$ \lap u = 0, \;\;\;\; x \in \mathbb{R}^n,$$
then it must be constant.}

One of its important applications is in the proof of the Fundamental Theorem of Algebra. It is also a key
ingredient in deriving point-wise a priori estimates for solutions to a family of elliptic equations in bounded domains with prescribed boundary values (see \cite{CLL1, GS, Z}), and thus obtaining the existence of solutions \cite {FLN, R}.

Liouville type theorem can also be used to study geometrical and reaction diffusion problems (see \cite{AW, CY, KN, KLS, M}), and to derive singularity and decay estimates (see \cite{PQS}).

A similar result has been established for $s$-harmonic functions \cite{BKN, Fa, CDL}:

{\em Assume that $0<s<1$. If $u$ is bounded from below or from above and
$$ (-\lap)^s u = 0, \;\;\;\; x \in \mathbb{R}^n,$$
then it must be constant.}

Here
$$(-\lap)^{s} u(x) = C_{n,s} \, \lim_{\epsilon \ra 0} \int_{R^n\setminus B_{\epsilon}(x)} \frac{u(x)-u(z)}{|x-z|^{n+2s}} dz.$$
And in order for the integral to converge, we require that $u \in C^{1,1}_{loc} \cap {\cal L}_{2s}$ with
$${\cal L}_{2s}=\{u \in L^1_{loc}  \mid \int_{R^n}\frac{|u(x)|}{1+|x|^{n+2s}} \, d x <\infty\}.$$

Besides the above mentioned applications, this Liouville theorem has also been used to prove the equivalences between fractional nonlinear equations and the corresponding integral equations (see \cite{CLM, ZCCY}), thus one can employ integral equations methods, such as {\em method of moving planes in integral forms} (\cite{CFY, CLO, CLO1, FC, HLZ, Lei, LLM, LZh}) to study qualitative properties of the solutions for the original fractional equations.

To study fractional harmonic functions, one powerful tool is the Poisson representation:

 $$
u(x)= \underset{|y|>r}{\int}P_{r}(y,x)u(y)dy, \;\; \mbox{ for } |x|<r,
$$
where $P_r(y,x)$ is the Poisson kernel:
$$
P_{r}(y,x)=\left\{\begin{array}{ll}
\frac{\Gamma(n/2)}{\pi^{\frac{n}{2}+1}} \sin\frac{\pi\alpha}{2}
\left[\frac{r^{2}-|x|^{2}}{|y|^{2}-r^{2}}\right]^{\frac{\alpha}{2}}\frac{1}{|x-y|^{n}},& |y|>r,\\
0,& |y|<r.
\end{array}
\right.
$$

If $u$ is bounded from one side, differentiating under the integral sign and letting $r \ra \infty$, one will be able to show that
$$\frac{\partial u}{\partial x_i} (x) = 0, \;\; i = 1, 2, \cdots, n, \; x \in \mathbb{R}^n,$$
and this leads to the conclusion that $u$ is constant (see \cite{BKN, CDL, Fa}). Actually more general results under weaker conditions were obtained in \cite{CDL} and \cite{Fa}.

The other important tool in studying fractional harmonic functions in $\mathbb{R}^n$ is the Fourier transform ${\cal F}$. From
$$ 0 = (-\lap)^s u(x) = {\cal F}((-\lap)^s u)(\xi) = |\xi|^{2s} {\cal F}(u)(\xi),$$
it follows immediately that
$$ {\cal F}(u)(\xi) = 0 \;\; \mbox{ for } \xi \neq 0.$$
 Hence ${\cal F} u$ is a finite combination of the Dirac's delta measure and its derivatives. Therefore $u$ is a polynomial. Under further restrictions that $u \in {\cal L}_{2s}$ and bounded from one side, it has to be constant.

Unfortunately, these two effective tools no longer work for the fractional $p$-Laplacians due to the full nonlinearity of the operator,
while this kind of operators have been recently used in many applications including continuum mechanics, population dynamics, and many different non-local diffusion
problems \cite{AMRT}. It is also applied to study the non-local ``Tug-of-War" game (see \cite{BCF, BCF1}).

As usual, one of the fundamental problem in studying the fractional $p$-Laplacians is the uniqueness of fractional $p$-harmonic functions. Since the conventional methods do not work, so far as we know, there has not been any results in this respect. In this paper, we approach such a problem in a completely different way. We prove that if $u$ is a bounded fractional $p$-harmonic functions, then it is symmetric about any given hyper-plane in $\mathbb{R}^n$, and therefore, it must be constant. We believe that this idea will be useful in studying other related problems.

To obtain the symmetry of $u$ with respect to a given hyper-plane $T$. Let $w(x)=u(\bar{x})-u(x)$ with $\bar{x}$ being the reflection of $x$ with respect to plane $T$. We prove that
$$ w(x) \leq 0 \;\;\; \mbox{ for all $x$ in the half space on one side of the plane }.$$
This actually is a {\em maximum principle } for anti-symmetric functions on a half space, an unbounded region, without assuming that the function vanishes near infinity.
  For convenience in applications, we summarize it as
\begin{mthm} (Maximum Principle)
Let $T$ be any given hyper-plane in $\mathbb{R}^n$ and $\Sigma$ be the half space on one side of the plane.  Let $$w(x)=u(\tilde{x})-u(x)$$ with $\tilde{x}$ being the reflection of $x$ with respect to plane $T$.

Assume that $w(x)$ is bounded in $\Sigma$ and
$$ (-\lap)^s_p u(\tilde{x}) - (\lap)^s_p u(x) \leq 0$$
at the points $x \in \Sigma$ where
$$u(\tilde{x}) > u(x).$$

Then
$$w(x) \leq 0, \;\;\;\; \forall \, x \in \Sigma.$$
\label{mthm2}
\end{mthm}

We believe that this {\em maximum principle} will become a powerful tool in carrying out the {\em method of moving planes} on unbounded domains. To illustrate this, let's study a simple example related to
\medskip

{\bf De Giorgi Conjecture} \cite{G}. If $u$ is a solution of
\be -\lap u(x) = u(x)-u^3(x), \;\;\; x \in \mathbb{R}^n
\label{DG1}
\ee
such that
$$|u(x)| \leq 1, \; \lim_{x_n \ra \pm \infty} u(x',x_n) = \pm 1 \mbox{ for all } x' \in \mathbb{R}^{n-1}, \mbox{ and } \frac{\partial u}{\partial x_n} >0.$$
Then there exists a vector ${\bf a} \in \mathbb{R}^{n-1}$ and a function $u_1 : \mathbb{R} \ra \mathbb{R}$ such that
$$ u(x', x_n) = u_1({\bf a} \cdot x' + x_n), \;\;\; \forall \, x \in \mathbb{R}^n.$$

Now we consider the fractional version of (\ref{DG1}):

\be
(-\lap)^s_p u(x) = f(u(x)), \;\;\;\;  x \in \mathbb{R}^n.
\label{DG}
\ee

Then applying Theorem \ref{mthm2}, we will be able to derive

\begin{mthm} Suppose that $u$ is a solution of (\ref{DG}), and
$$|u(x)| \leq 1, \;\;\; \forall \, x \in \mathbb{R}^n, $$
\be \lim_{x_n \ra \pm \infty} u(x',x_n) =  \pm 1 \mbox{ uniformly for all } x' \in \mathbb{R}^{n-1},
\label{u1}
\ee
and
\be
f(t) \mbox{ is non-increasing for $|t|$ sufficiently close to $1$},
\label{f}
\ee
then there exists $M>0$, such that
\be
\frac{\partial u}{\partial x_n} > 0, \;\;\; \mbox{ forall  $x$ with } |x_n| \geq M.
\label{u2}
\ee
\label{mthm3}
\end{mthm}

\begin{mrem}

Note here the condition on $f$ is satisfied for $f(u)=u-u^3$ as given in the {\em De Giorgi Conjecture}.
However, our condition (\ref{u1}) is stronger than that in the {\em conjecture}.

\end{mrem}

To see how Theorem \ref{mthm2} implies Theorem \ref{mthm3}, let
$$T_{\lambda} = \{ x \in \mathbb{R}^n \mid x_n = \lambda \}$$ and
$$\Sigma_\lambda = \{ x \in \mathbb{R}^n \mid x_n > \lambda \}. $$
For $x \in \Sigma_\lambda$, let $x^\lambda$ be the reflection of $x$ with respect to plane $T_\lambda$.
Let $$w_\lambda (x) = u(x^\lambda)-u(x).$$
Then for $\lambda$ sufficiently large, (\ref{u1}) and (\ref{f}) imply that
$$ (-\lap)^s_p u(x^\lambda) - (-\lap)^s_p u(x) \leq 0$$
at the points $x \in \Sigma_\lambda$ where
$$u(x^\lambda) > u(x).$$
Now it follows from Theorem \ref{mthm2} that

\be
w_\lambda (x) \leq 0, \;\;\;\; \forall \, x \in \Sigma_\lambda, \mbox{ for all sufficiently large } \lambda.
\label{w1}
\ee
Consequently, a standard arguments will lead to (\ref{u2}) for $x_n \geq M.$ Similarly, one can show (\ref{u2}) for $x_n \leq -M$.

Inequality (\ref{w1}) actually provides a starting point to move the plane $T_\lambda$. If we can move the plane all the way down, then we prove that
$$\frac{\partial u}{\partial x_n} > 0, \;\;\; \forall \, x \in \mathbb{R}^n.$$
For this purpose, we need to establish a {\em narrow region principle} on unbounded domains (without assuming that the function vanishes near infinity). This requires further condition on $f(\cdot)$, and fortunately, it is satisfied by the model function $f(u)=u-u^3$.  We will present this {\em narrow region principle} in our next paper. Notice that, previously in carrying out the {\em methods of moving planes on unbounded regions}, people usually assume that the function $u$ vanishes near infinity, or make a Kelvin transform so that the transformed function vanishes near infinity, or consider $u(x)/g(x)$ for a proper choice of $g(x)$, so that $u(x)/g(x)$ vanishes near infinity. Unfortunately, the later two approaches do not work on the nonlinear nonlocal operators such as the fractional $p$-Laplacian. Hence it is useful to develop a {\em method of moving planes} for such operators that applies to unbounded regions while only assuming the function be bounded. As illustrated above, our Theorem \ref{mthm2} provides a starting point to move the planes in such a situation. Therefore, we believe that the result, and in particular, the ideas in its proof, will become useful tools in studying many other problems involving nonlinear and nonlocal operators on unbounded domains.

In Sections 2, we will establish the {\em maximum principle}  and hence prove Theorem \ref{mthm2}. In Section 3, we will use the {\em maximum principle} to derive the Liouville Theorem (Theorem \ref{main}) and Theorem \ref{mthm3}. Two preliminary lemmas will be proved in Section 4.

For more articles concerning Liouville theorems, please see \cite{PQS, LZ, CFY, MM} and the references therein.

\section{The proof of the {\em maximum principle}}

In this section, we prove

\begin{thm} (Maximum Principle)
Let $T$ be any given hyper-plane in $\mathbb{R}^n$ and $\Sigma$ be a region to one side of the plane.  Let $w(x)=u(\tilde{x})-u(x)$ with $\tilde{x}$ being the reflection of $x$ with respect to plane $T$.

Assume that
\be
(-\lap)^s_p u(\tilde{x}) - (-\lap)^s_p u(x) \leq 0
\label{a10}
\ee
at the points $x \in \Sigma$ where
$$u(\tilde{x}) > u(x).$$

Then
\be
w(x) \leq 0, \;\;\;\; \forall \, x \in \Sigma.
\label{w}
\ee
\label{thm2}
\end{thm}

{\em Outline of the proof.}

 Suppose (\ref{w}) is false, we have
\begin{eqnarray}\label{A}
A:= \sup_{\Sigma}w (x)>0.
\end{eqnarray}

If this supremum can be attained, say at point $x^o$, then one can derive
immediately
$$ (-\lap)^s_p u (\tilde{x^o}) - (-\lap)^s_p u(x^o) > 0$$
(see \cite{CL}, the proof of {\em maximum principles for anti-symmetric functions}).
This would contradict (\ref{a10}).

Now the main problem is that $\Sigma$ is an unbounded region, and hence
$\mathop{\sup}\limits_{\Sigma}w_ (x)$ may not be attained. To circumvent
this difficulty, we choose a point $x^o$ at which $w$ is very close to $A$. Select
a compactly supported radially symmetric function $\psi (x) = \psi(|x-x^o|)$ which is strictly monotone decreasing with respect to $|x-x^o|$ so that the function
$$ w (x) + \epsilon (\psi (\tilde{x}) - \psi (x))$$
attains its maximum at some point $\bar{x}$ near $x^o$. Here
$\psi (\tilde{x}) - \psi (x)$ is an anti-symmetric function. Then we will be able to estimate
$$ Q := (-\lap)^s_p ( u + \epsilon \psi )(\tilde{x}) - (-\lap)^s_p ( u + \epsilon \psi )(x)$$
at the maximum point $\bar{x}$.

On one hand, by using the integral definition of $(-\lap)^s_p$ and the anti-symmetry property of $w (x) + \epsilon (\psi (\tilde{x}) - \psi (x))$, we derive
\be
Q \geq C_1 A^{p-1}
\label{a2}
\ee
for some positive constant $C_1$.

On the other hand,
from (\ref{a10}) and the bounded-ness of $(-\lap)^s_p \psi$, we will show that
$$
Q \mbox{ can be made as small as we wish provided  $\epsilon$ is sufficiently small.}
$$
This will contradict (\ref{a2}) since here $A$ is a fixed positive number, while $\epsilon$
can be chosen as small as we wish provided $w (x^o)$ is sufficiently close to $A$.
\medskip

{\em Now we carry out the details of the proof.}
\smallskip

For any $0<\gamma<1,$ there exists $x_0 \in \Sigma$, such that
$$
w (x_0)\geq \gamma A.
$$

By translation, rotation, and re-scaling, we may assume that
$$ T = \{ x \mid x_1=0 \},  \;\; \Sigma =\{x: \mid x_1<0 \}, \; \mbox{ and } x_0=(-2,0, \cdots ,0).$$
Then
$$\tilde{x} = (-x_1, x_2, \cdots, x_n).$$
Denote
$$u_0(x) = u(\tilde{x}).$$
Then
$$
w(x)=u_0(x)-u(x).
$$
Let
$$
\eta(x)=\left\{\begin{array}{ll}
a e^{\frac{1}{|x|^2-1}}, &|x|<1 \\
  0,& |x|\geq 1,
\end{array} \right.
$$
taking $a=e$ such that $\eta (0)= \max_{\mathbb{R}^n}\eta (x)=1.$ It is well-known that this $\eta(x)$ is a $C_0^\infty$ function, and hence
$$|(-\lap)^s_p \eta(x)| \leq C, \;\; x \in \mathbb{R}^n.$$
Also it is monotone decreasing with respect to $|x|$.

Set
$$
\psi (x) =   \eta (x+x_0) \mbox{ and } \psi_0(x)=\psi(\tilde{x})=\eta (x-x_0).
$$
Then $\psi_0(x)-\psi(x)$ is an anti-symmetric function with respect to the plane $T$.

Now pick $\varepsilon>0$ such that
$$
w(x_0) +\varepsilon \psi_0(x_0)-\varepsilon  \psi(x_0)\geq A,
$$
notice that $\psi(x_0)=0$ and $\psi_0(x_0)=1$, hence we only need to let
\begin{eqnarray}\label{varepsilon}
 \varepsilon=(1-\gamma)A.
\end{eqnarray}
Since for any $x\in \Sigma \backslash B_1(x_0),$ $w(x)\leq A$ and $ \psi_0(x)=\psi (x)=0,$ we have
$$
w(x_0) +\varepsilon \psi_0(x_0)-\varepsilon  \psi(x_0)\geq w(x) +\varepsilon \psi_0(x)-\varepsilon \psi(x), \text{~for~any~} x \in \Sigma \backslash B_1(x_0),
$$
which means that the supremum of the function $w(x) +\varepsilon \psi_0(x)-\varepsilon \psi(x)$ in $\Sigma$ is achieved on $\overline{B_1(x_0)},$ hence there exists a point $\bar{x} \in \overline{B_1(x_0)}$ such that
\begin{eqnarray}\label{wmaximum}
 w(\bar{x}) +\varepsilon \psi_0(\bar{x})-\varepsilon \psi(\bar{x})=\max_{\Sigma}( w(x) +\varepsilon \psi_0(x)-\varepsilon \psi(x))\geq A.
\end{eqnarray}

Now we will evaluate the lower bound and the upper bound of
\begin{eqnarray}\label{estimate}
 (-\Delta)_p^s(u_0+\varepsilon \psi_0)(\bar{x})-(-\Delta)_p^s(u+\varepsilon \psi)(\bar{x})
\end{eqnarray}
respectively to derive a contradiction.

We first estimate the lower bound of (\ref{estimate}) by direct calculations.

For convenience, we denote $G(t)=|t|^{p-2}t,$ then $G'(t)=(p-1)|t|^{p-2}\geq 0$
and may assume that $C_{n, sp}=1$. We have

\begin{eqnarray}\label{calculation}
&&(-\Delta)_p^s(u_0+\varepsilon \psi_0)(\bar{x})-(-\Delta)_p^s(u+\varepsilon \psi)(\bar{x})\nonumber\\
&=&PV\int_{\mathbb{R}^n}\frac{1}{|\bar{x}-y|^{n+sp}}\left[G(u_0(\bar{x})+\varepsilon \psi_0(\bar{x})-u_0(y)-\varepsilon \psi_0(y))\right.\nonumber\\
&&~ \quad\quad\quad\quad\quad\quad\quad\quad\left.-G(u(\bar{x})+\varepsilon \psi (\bar{x})-u(y)-\varepsilon \psi(y))\right]dy\nonumber\\
&=&PV\int_{\Sigma}\frac{1}{|\bar{x}-y|^{n+sp}}\left[G(u_0(\bar{x})+\varepsilon \psi_0(\bar{x})-u_0(y)-\varepsilon \psi_0(y))\right. \nonumber\\
&&~\quad\quad\quad\quad\quad\quad \quad\quad \left.-G(u(\bar{x})+\varepsilon \psi (\bar{x})-u(y)-\varepsilon \psi(y))\right]dy\nonumber\\
&&+PV\int_{\Sigma}\frac{1}{|\bar{x}+y|^{n+sp}}\left[G(u_0(\bar{x})+\varepsilon \psi_0(\bar{x})-u(y)-\varepsilon \psi(y))\right. \nonumber\\
&&\quad\quad\quad\quad\quad\quad\quad\quad~~~\left.-G(u(\bar{x})+\varepsilon \psi (\bar{x})-u_0(y)-\varepsilon \psi_0(y))\right]dy\nonumber\\
&=&PV\int_{\Sigma}\left(\frac{1}{|\bar{x}-y|^{n+sp}}-\frac{1}{|\bar{x}+y|^{n+sp}}\right)\left[G(u_0(\bar{x})+\varepsilon \psi_0(\bar{x})-u_0(y)-\varepsilon \psi_0(y))\right.\nonumber\\
&&\quad\quad\quad\quad\quad\quad\quad\quad\quad\quad\quad\quad\quad\quad\quad\quad\left.-G(u(\bar{x})+\varepsilon \psi (\bar{x})-u(y)-\varepsilon \psi(y))\right]dy\nonumber\\
 &&+PV\int_{\Sigma}\frac{1}{|\bar{x}+y|^{n+sp}}
\left\{ \left[G(u_0(\bar{x})+\varepsilon \psi_0(\bar{x})-u(y)-\varepsilon \psi(y))\right.\right.\nonumber\\
 &&\quad\quad\quad\quad\quad\quad\quad\quad\quad~~
 \left.\left. -G(u(\bar{x})+\varepsilon \psi (\bar{x})-u_0(y)-\varepsilon \psi_0(y))\right]\right.\nonumber\\
 &&\quad\quad\quad\quad\quad\quad\quad\quad\quad~~\left.+\left[G(u_0(\bar{x})+\varepsilon \psi_0(\bar{x})-u_0(y)-\varepsilon \psi_0(y))\right.\right.\nonumber\\
 &&\quad\quad\quad\quad\quad\quad\quad\quad\quad~~
 \left.\left.-G(u(\bar{x})+\varepsilon \psi (\bar{x})-u(y)-\varepsilon \psi(y))\right]\right\}dy\nonumber\\
 &=&PV\int_{\Sigma}\left(\frac{1}{|\bar{x}-y|^{n+sp}}-\frac{1}{|\bar{x}+y|^{n+sp}}\right)\left[G(u_0(\bar{x})+\varepsilon \psi_0(\bar{x})-u_0(y)-\varepsilon \psi_0(y))\right.\nonumber\\
&&\quad\quad\quad\quad\quad\quad\quad\quad\quad\quad\quad\quad\quad\quad\quad\quad \left.-G(u(\bar{x})+\varepsilon \psi (\bar{x})-u(y)-\varepsilon \psi(y))\right]dy\nonumber\\
 &&+PV\int_{\Sigma}\frac{1}{|\bar{x}+y|^{n+sp}}\left\{\left[G(u_0(\bar{x})+\varepsilon \psi_0(\bar{x})-u(y)-\varepsilon \psi(y)) \right.\right.\nonumber\\
 &&\quad\quad\quad\quad\quad\quad\quad\quad\quad\quad
 \left.\left.-G(u(\bar{x})+\varepsilon \psi (\bar{x})-u(y)-\varepsilon \psi(y))
 \right]\right.\nonumber\\
 &&\quad\quad\quad\quad\quad\quad\quad\quad\quad\quad\left.+\left[G(u_0(\bar{x})+\varepsilon \psi_0(\bar{x})-u_0(y)-\varepsilon \psi_0(y))\right.\right.\nonumber\\
 &&\quad\quad\quad\quad\quad\quad\quad\quad\quad\quad\left.\left.-G(u(\bar{x})+\varepsilon \psi (\bar{x})-u_0(y)-\varepsilon \psi_0(y))\right]\right\}dy\nonumber\\
 &:=&I_1+I_2.
\end{eqnarray}

To estimate $I_1$, we first notice that
$$
\frac{1}{|\bar{x}-y|^{n+sp}}-\frac{1}{|\bar{x}+y|^{n+sp}} >0, \text{~for~all~} y \in \Sigma,
$$
while for the second part in the integral, for any $y \in \Sigma,$ we have
\begin{eqnarray*}
G(u_0(\bar{x})+\varepsilon \psi_0(\bar{x})-u_0(y)-\varepsilon \psi_0(y))-G(u(\bar{x})+\varepsilon \psi (\bar{x})-u(y)-\varepsilon \psi(y))\geq0
\end{eqnarray*}
due to the strict monotonicity of $G$ and the fact that
\begin{eqnarray*}
&&(u_0(\bar{x})+\varepsilon \psi_0(\bar{x})-u_0(y)-\varepsilon \psi_0(y))-(u(\bar{x})+\varepsilon \psi (\bar{x})-u(y)-\varepsilon \psi(y))\\
 &=&(w(\bar{x})+\varepsilon \psi_0(\bar{x})-\varepsilon \psi(\bar{x}))-(w(y)+\varepsilon \psi_0(y)-\varepsilon \psi(y))\\
 &\geq&0.
\end{eqnarray*}
Therefore,
\begin{align}\label{I1}
I_1&=\int_{\Sigma}\left(\frac{1}{|\bar{x}-y|^{n+sp}}-\frac{1}{|\bar{x}+y|^{n+sp}}\right)\left[G(u_0(\bar{x})+\varepsilon \psi_0(\bar{x})-u_0(y)-\varepsilon \psi_0(y))\right.\nonumber\\
&\quad\quad\quad \quad\quad\quad \quad\quad\quad\quad \quad\quad\quad \quad\quad\left.-G(u(\bar{x})+\varepsilon \psi (\bar{x})-u(y)-\varepsilon \psi(y))\right]dy\nonumber\\
&\geq0.
\end{align}

To estimate the integral $I_2$ in (\ref{calculation}), we need the following analysis lemma.
\begin{mlem}\label{LeG0}
For $G(t)=|t|^{p-2}t,$  there exists a constant $C>0$ such that
\begin{eqnarray}\label{eqG0}
G(t_2)-G(t_1) \geq C(t_2 -t_1)^{p-1}
\end{eqnarray}
for any $t_2>t_1$.
\end{mlem}

The proof of this lemma is quite elementary, while for reader's convenience, we include it in Section 4.

Now we apply this lemma to estimate $I_2$. By (\ref{wmaximum}), for any $y \in \Sigma,$ we have
\begin{eqnarray*}
&&(u_0(\bar{x})+\varepsilon \psi_0(\bar{x})-u(y)-\varepsilon \psi(y))-
(u(\bar{x})+\varepsilon \psi (\bar{x})-u(y)-\varepsilon \psi(y))\nonumber\\
&=&w(\bar{x})+\varepsilon \psi_0(\bar{x})-\varepsilon \psi(\bar{x})\geq A,
\end{eqnarray*}
and
\begin{eqnarray*}
 &&(u_0(\bar{x})+\varepsilon \psi_0(\bar{x})-u_0(y)-\varepsilon \psi_0(y))-(u(\bar{x})+\varepsilon \psi (\bar{x})-u_0(y)-\varepsilon \psi_0(y))\nonumber\\
 &=&w(\bar{x})+\varepsilon \psi_0(\bar{x})-\varepsilon \psi(\bar{x})\geq A.
\end{eqnarray*}
Then by Lemma \ref{LeG0}, we derive
\begin{eqnarray}\label{I2}
 I_2&=&\int_{\Sigma}\frac{1}{|\bar{x}+y|^{n+sp}}\left\{\left[G(u_0(\bar{x})+\varepsilon \psi_0(\bar{x})-u(y)-\varepsilon \psi(y)) \right.\right.\nonumber\\
 &&\quad\quad\quad\quad\quad\quad\quad
 \left.\left.-G(u(\bar{x})+\varepsilon \psi (\bar{x})-u(y)-\varepsilon \psi(y))
 \right]\right.\nonumber\\
 &&\quad\quad\quad\quad\quad\quad\quad\left.+\left[G(u_0(\bar{x})+\varepsilon \psi_0(\bar{x})-u_0(y)-\varepsilon \psi_0(y))\right.\right.\nonumber\\
 &&\quad\quad\quad\quad\quad\quad\quad
 \left.\left.-G(u(\bar{x})+\varepsilon \psi (\bar{x})-u_0(y)-\varepsilon \psi_0(y))\right]\right\}dy\nonumber\\
 &\geq& C \int_{\Sigma}\frac{1}{|\bar{x}+y|^{n+sp}}(w(\bar{x})+\varepsilon \psi_0(\bar{x})-\varepsilon \psi(\bar{x}))^{p-1} dy\nonumber\\
 &\geq& C \int_{\Sigma}\frac{1}{|\bar{x}+y|^{n+sp}} A^{p-1}dy\nonumber\\
 &\geq& C_1 A^{p-1}.
\end{eqnarray}

 Combining (\ref{calculation}), (\ref{I1}) and (\ref{I2}), we deduce
\begin{eqnarray}\label{geq}
(-\Delta)_p^s(u_0+\varepsilon \psi_0)(\bar{x})-(-\Delta)_p^s(u+\varepsilon \psi)(\bar{x})\geq C_1 A^{p-1}.
\end{eqnarray}

To derive a contradiction, we also estimate the upper bound of (\ref{estimate}). To proceed, we need the following lemma, which will be proved in Section 4.

 \begin{mlem}\label{Convergence}
 Assume that  $u \in C_{loc}^{1,1} \cap {\cal L}_{sp}(\mathbb{R}^n)$ and $\psi \in C_0^\infty(\mathbb{R}^n),$ then for all small $\delta >0$, we have
\begin{eqnarray*}
\mid(-\Delta)_p^s(u+\varepsilon \psi)-(-\Delta)_p^su \mid \leq \varepsilon C_\delta +C \delta^{p(1-s)},
\end{eqnarray*}
where $C$ is independent of $\varepsilon$ while $C_\delta$ may depend on $\delta$.
\end{mlem}

Now we employ this lemma to continue evaluating the upper bound of (\ref{estimate}), we derive
\begin{eqnarray}\label{upper bound}
 &&(-\Delta)_p^s(u_0+\varepsilon \psi_0)(\bar{x})-(-\Delta)_p^s(u+\varepsilon \psi)(\bar{x})\nonumber\\
 &\leq& (-\Delta)_p^s u_0(\bar{x})-(-\Delta)_p^s u(\bar{x})+ 2 \varepsilon C_\delta +C_2 \delta^{p(1-s)}\nonumber\\
 &\leq & 2 \varepsilon C_\delta +C_2 \delta^{p(1-s)}.
\end{eqnarray}
Combining (\ref{geq}) and (\ref{upper bound}), we derive
$$C_1 A^{p-1} \leq \varepsilon c_\delta +C_2 \delta^{p(1-s)}.$$
 We first choose $\delta$ small such that
 $$
 C_2 \delta^{p(1-s)}\leq \frac{C_1}{3}A^{p-1},
 $$
then for such $\delta$, let $\gamma$ be sufficiently close to $1$, hence $\varepsilon =(1-\gamma)A$ is small such that
$$
\varepsilon c_\delta \leq \frac{C_1}{3}A^{p-1},
$$
which contradicts with $A>0.$
This implies (\ref{w}) and thus completes the proof of the theorem.

\section{Applications of the {\em maximum principle}--the Liouville theorem and more}

 In this section, we will use the {\em maximum principle} established in the previous section to prove Theorem \ref{main} and \ref{mthm3}.
 \medskip

{\bf The proof of Theorem \ref{main}}
\medskip

We show that $u$ is symmetric with respect to any hyper plane. To this end, let $x_n$ be any given direction in $\mathbb{R}^n$ and let
 $$
T_\lambda=\{x \in \mathbb{R}^n \mid x_n=\lambda \text{~for~}\lambda \in \mathbb{R}\}
$$
be a plane perpendicular to $x_n$-axis.

Let
$$
\Sigma_\lambda=\{x \in \mathbb{R}^n \mid x_n > \lambda\}
$$
be the region above the plane $T_\lambda$.

For $x \in \Sigma_\lambda$, let
$$
x^\lambda=(x_1, x_2, \cdots, 2 \lambda- x_n)
$$
be its reflection about the plane $T_\lambda.$

Denote
$$
w_\lambda(x)=u(x^\lambda)-u(x).
$$

Applying Theorem \ref{mthm2}, we arrive immediately
$$
 w_\lambda(x) \leq 0 \text{~in~} \Sigma_\lambda,
$$

Similarly, we can prove that
$$
 w_\lambda(x) \geq 0 \text{~in~} \Sigma_\lambda,
$$
therefore,
\begin{eqnarray}
 w_\lambda(x) \equiv 0 \text{~in~} \Sigma_\lambda.
 \label{w=0}
\end{eqnarray}
These imply that  $u(x)$ is symmetric with respect to plane $T_\lambda$ for any $\lambda \in \mathbb{R}$.

Since the $x_n$-direction can be chosen arbitrarily, (\ref{w=0}) implies $u$ is radially symmetric about any point, it follows that
$$
u(x) \equiv C.
$$

This completes the proof of the theorem.
\medskip

{\bf The proof of Theorem \ref{mthm3}}.
\smallskip

We only need to prove that
$$
w_\lambda (x)=u(x^\lambda) -u(x)\leq 0 \mbox{ in } \Sigma_\lambda, \mbox{ for sufficiently large } \lambda.
$$
If not, then
\be
\sup_{\Sigma_\lambda} w_\lambda(x)=A>0.
\label{supw}
\ee
Therefore, for any $0<\gamma<1$, there exists $x_0 \in \Sigma_\lambda$ such that
$$
w(x_0)\geq \gamma A.
$$
 By re-scaling, we may assume that $dist\{x_0, T_\lambda\}=2.$ Let
$$
\eta(x)=\left\{\begin{array}{ll}
a e^{\frac{1}{|x|^2-1}}, &|x|<1 \\
  0,& |x|\geq 1.
\end{array} \right.
$$
Let $a=e$ such that $\eta(0)= \mathop{\max}\limits_{\mathbb{R}^n}\eta(x)=1.$
Set
$$
\psi (x) =  \eta (x-x^\lambda_0) \mbox{ and } \psi_\lambda (x)=\eta (x-x_0).
$$
Then $\psi_\lambda(x)-\psi(x)$ is an anti-symmetric function with respect to the plane $T_\lambda$.

Now pick $\varepsilon=(1-\gamma)A>0$ such that
$$
w_\lambda(x_0) + \varepsilon \psi_\lambda(x_0)-\varepsilon  \psi(x_0)\geq A.
$$
It follows that there exists a point $\bar{x} \in \overline{B_1(x_0)}$ such that
\begin{eqnarray*}
 w_\lambda(\bar{x}) + \varepsilon \psi_\lambda(\bar{x}) - \varepsilon \psi(\bar{x})=\max_{\Sigma}( w_\lambda(x) + \varepsilon \psi_\lambda(x) - \varepsilon \psi(x))\geq A.
\end{eqnarray*}

Similar to the proof of Theorem \ref{thm2}, we will be able to estimate
$$
Q:=(-\lap)_p^s(u_\lambda +\varepsilon \psi_\lambda)(\bar{x})-(-\lap)_p^s(u + \varepsilon \psi)(\bar{x})
$$
at the maximum point $\bar{x}.$

On one hand, we have
$$
Q\geq CA^{p-1}.
$$

On the other hand, since
$$
 w_\lambda(\bar{x}) +\varepsilon \psi_\lambda(\bar{x})-\varepsilon \psi(\bar{x})\geq w_\lambda(x_0) +\varepsilon \psi_\lambda(x_0)-\varepsilon  \psi(x_0),
$$
 $\psi(\bar{x})=\psi (x_0)=0$ and $\psi_\lambda(x_0) \geq \psi_\lambda (\bar{x})$,
we have
$$w_\lambda (\bar{x})\geq w_\lambda (x_0)>0.$$
Hence,
$$
u_\lambda (\bar{x}) >u (\bar{x}).
$$
By the monotonicity of $f$ and Lemma \ref{Convergence}, we derive
\begin{eqnarray*}
 Q &\leq& (-\lap)_p^s u_\lambda (\bar{x})- (-\lap)_p^s u (\bar{x})+ \varepsilon c_\delta +C_2 \delta^{p(1-s)}\\
 &\leq&f(u_\lambda (\bar{x}))-f(u (\bar{x}))+ \varepsilon c_\delta +C_2 \delta^{p(1-s)}\\
 &\leq& \epsilon c_\delta +C_2 \delta^{p(1-s)}.
\end{eqnarray*}
Therefore, we derive a contradiction.

This completes the proof of Theorem \ref{mthm3}.

\section{Preliminary Lemmas}

\begin{mlem}\label{LeG}
For $G(t)=|t|^{p-2}t, p>2,$  there exists a constant $C>0$ such that
\begin{eqnarray}\label{eqG}
G(t_2)-G(t_1) \geq C(t_2 -t_1)^{p-1}
\end{eqnarray}
for any $t_2\geq t_1$.
\end{mlem}

\textbf{Proof.} We consider three possible cases.

Case i) \,\, $t_2 \geq t_1\geq 0.$

If $t_1 \geq \frac{t_2}{2},$ then
\begin{eqnarray*}
G(t_2)-G(t_1)&=&t_2^{p-1}-t_1^{p-1}\\
&=&(p-1)\xi^{p-2}(t_2-t_1)\\
&\geq& (p-1)\left(\frac{t_2}{2}\right)^{p-2}(t_2-t_1)\\
&\geq& C(t_2 -t_1)^{p-1},
\end{eqnarray*}
where $\xi$ lies between $t_1$ and $t_2,$ this implies (\ref{eqG}).

If $0<t_1 \leq \frac{t_2}{2},$ then
\begin{eqnarray*}
G(t_2)-G(t_1)=t_2^{p-1}-t_1^{p-1}\geq t_2^{p-1} -\left(\frac{t_2}{2}\right)^{p-1}\geq C(t_2 -t_1)^{p-1}.
\end{eqnarray*}
This implies (\ref{eqG}).

Case ii)  $t_1\leq t_2 \leq 0.$

In this case, $|t_1| \geq |t_2|\geq 0$ and we can derive from Case i) that
\begin{eqnarray*}
G(t_2)-G(t_1)=|t_1|^{p-1}-|t_2|^{p-1}\geq C(|t_1| -|t_2|)^{p-1}=C(t_2 -t_1)^{p-1}.
\end{eqnarray*}
This implies (\ref{eqG}).

Case iii) \,\, $t_1\leq 0 \leq t_2.$

In this case, $t_1$ and $t_2$ are of the different signs, then
\begin{eqnarray*}
G(t_2)-G(t_1)=|t_2|^{p-1}+|t_1|^{p-1}\geq C(|t_2| +|t_1|)^{p-1}=C(t_2 -t_1)^{p-1}.
\end{eqnarray*}
This implies (\ref{eqG}), and hence completes the proof of the lemma.

 \begin{mlem}\label{convergence}
 Assume that  $u \in C_{loc}^{1,1} \cap L_{sp}(\mathbb{R}^n)$ and $\psi \in C_0^\infty(\mathbb{R}^n),$ then for all small $\delta >0$, we have
\begin{eqnarray}\label{Iconvergence}
\mid (-\Delta)_p^s(u+\varepsilon \psi)-(-\Delta)_p^su \mid \leq \varepsilon C_\delta +C \delta^{p(1-s)},
\end{eqnarray}
where $C$ is independent of $\varepsilon$ while $C_\delta$ may depend on $\delta$.
\end{mlem}
\textbf{Proof.}
Denote $v_\varepsilon (x)=u(x)+\varepsilon \psi (x),$ then
\begin{eqnarray*}
&&(-\Delta)_p^s(u+\varepsilon \psi)(x)-(-\Delta)_p^su(x)\\
&=&PV \int_{\mathbb{R}^n}\frac{G(v_\varepsilon(x)-v_\varepsilon(y))}{|x-y|^{n+sp}}dy-PV\int_{\mathbb{R}^n}\frac{G(u(x)-
u(y))}{|x-y|^{n+sp}}dy,
\end{eqnarray*}
we divide $\mathbb{R}^n$ into two regions: $B_\delta(x)$ and $B_\delta^c(x)$.

(i) In $B_\delta^c(x)$.

 Since
\begin{eqnarray*}
&&G(v_\varepsilon(x)-v_\varepsilon(y))-G(u(x)-u(y))\nonumber\\
&=&\varepsilon (p-1)(\psi(x)-\psi(y))\int_0^1|u(x)-u(y)+t\varepsilon (\psi(x)-\psi(y))|^{p-2}dt\nonumber\\
&:=&\varepsilon F(x,y),
\end{eqnarray*}
and $u \in C_{loc}^{1,1} \cap L_{sp}(\mathbb{R}^n)$ and $\psi \in C_0^\infty(\mathbb{R}^n),$ we have
\begin{eqnarray}\label{EXTERIOR}
 \left|PV \int_{B_\delta^c(x)}
 \frac{\varepsilon F(x,y)}
 {|x-y|^{n+sp}}dy \right|\leq \varepsilon C_\delta.
\end{eqnarray}

(ii) In the ball $B_\delta(x)$.

In this case, we need the following basic inequality
\begin{eqnarray}\label{EQUALITY}
\mid|a+b|^{p-2}(a+b)-|a|^{p-2}a\mid \leq C (|a|+|b|)^{p-2}|b|,
\end{eqnarray}
which can be easily derived from the {\em mean value theorem}.

First, by Taylor expansion, we have
$$
v_\varepsilon(x)-v_\varepsilon(y)=\nabla v_\varepsilon(x)\cdot(x-y) +O(|x-y|^2).
$$
The anti-symmetry of $\nabla v_\varepsilon(x)\cdot(x-y)$ for  $y \in B_\delta(x)$ implies that
$$
PV \int_{B_\delta(x)}\frac{G(\nabla v_\varepsilon(x)\cdot(x-y))} {|x-y|^{n+sp}}dy=0.
$$
By (\ref{EQUALITY}), we obtain
\begin{eqnarray}\label{v-varepsilon}
&&\left|PV \int_{B_\delta(x)}\frac{G(v_\varepsilon(x)-v_\varepsilon(y))}{|x-y|^{n+sp}}dy
\right|\nonumber\\
&=&\left|PV \int_{B_\delta(x)}\frac{G(v_\varepsilon(x)-v_\varepsilon(y))}{|x-y|^{n+sp}}dy
-PV \int_{B_\delta(x)}\frac{G(\nabla v_\varepsilon(x)\cdot(x-y))}{|x-y|^{n+sp}}dy\right|\nonumber\\
&\leq & C \int_{B_\delta(x)}\frac{(\mid \nabla v_\varepsilon(x)\cdot(x-y)\mid + O(|x-y|^2) )^{p-2}O(|x-y|^2)}{|x-y|^{n+sp}}dy\nonumber\\
&=& C \int_{B_\delta(x)}\frac{ O(|x-y|^p) }{|x-y|^{n+sp}}dy\nonumber\\
&\leq& C_1 \delta^{p(1-s)},
\end{eqnarray}
where $C_1$ is independent of $\varepsilon,$ since for any fixed $x,$ we have
$$
|\nabla v_\varepsilon(x)|\leq |\nabla u(x)|+\varepsilon |\nabla \psi(x)|\leq |\nabla u(x)|+ |\nabla \psi(x)|, \text{~if~} \varepsilon \leq 1.
$$
Similarly, we have
\begin{eqnarray}\label{u}
\left|PV \int_{B_\delta(x)}\frac{G(u(x)-u(y))}{|x-y|^{n+sp}}dy\right|\leq C_2 \delta^{p(1-s)}.
\end{eqnarray}
Combining (\ref{EXTERIOR}), (\ref{v-varepsilon}) and (\ref{u}), we derive
$$|(-\Delta)_p^s(u+\varepsilon \psi)-(-\Delta)_p^su|\leq \varepsilon C_\delta +C_3 \delta^{p(1-s)}.$$

This completes the proof of Lemma \ref{convergence}.

\bigskip

{\em Authors' Addresses and E-mails:}
\medskip

Wenxiong Chen

Department of Mathematical Sciences

Yeshiva University

New York, NY, 10033 USA

wchen@yu.edu
\medskip

Leyun Wu

School of Mathematical Sciences and Institute of Natural Sciences,

Shanghai Jiao Tong University

Shanghai, China and

Department of Mathematical Sciences

Yeshiva University

New York, NY, 10033 USA

leyunwu@mail.nwpu.edu.cn
\medskip

\end{document}